\documentclass[10pt,reqno]{amsart}
\usepackage{amsmath,amssymb,amsfonts,amsthm,mathrsfs,enumitem}
\allowdisplaybreaks

\newtheorem{thm}{Theorem}[section]

\newtheorem{cor}[thm]{Corollary}
\newtheorem{lem}[thm]{Lemma}
\newtheorem{prop}[thm]{Proposition}

\theoremstyle{definition}
\newtheorem{defn}[thm]{Definition}
\theoremstyle{remark}
\newtheorem{rem}[thm]{Remark}
\newtheorem{exa}[thm]{Example}

\numberwithin{equation}{section}

\newcommand{\Z}{\mathbb Z}
\newcommand{\eps}{\varepsilon}
\newcommand{\supp}{\operatorname{supp}}
\newcommand{\dens}{\operatorname{dens}}
\newcommand{\Sub}{\operatorname{Sub}}
\newcommand{\Law}{\mathcal L}
\newcommand{\rev}[1]{#1^{\mathrm{rev}}}
\newcommand{\one}{\mathbf 1}

\newcommand{\bits}{\{0,1\}^{*}}

\begin{document}

\title[Cusick's conjecture]
{A first-exit proof of Cusick's sum-of-digits conjecture}
\author{Kaimin Cheng}
\address{School of Mathematical Sciences, China West Normal University, Nanchong 637002, P. R. China}
\email{ckm20@126.com}
\subjclass[2020]{Primary 11A63; Secondary 05A05, 05A20, 60G40}
\keywords{Binary sum of digits, Cusick conjecture, first-exit medians, subsequence ideals, deletion-closed languages, stopped random walks, digit reversal}
\date{}

\begin{abstract}
We prove Cusick's conjecture on the binary sum-of-digits function.  More
precisely, for every integer $t\ge1$ we show that
\[
c_t:=\lim_{N\to\infty}\frac1N
\#\{0\le n<N:\ s_2(n+t)\ge s_2(n)\}
>
\frac12,
\]
and in fact obtain the explicit bound
\[
c_t\ge \frac12+2^{-2s_2(t)-1},
\]
where $s_2(m)$ denotes the number of ones in the binary expansion of $m$. The proof is based on an exact deconvolution which replaces the distribution
of $s_2(n+t)-s_2(n)$ by a finite stopped random-walk law.  The required bias
is then proved through first-exit medians for principal subsequence ideals.
\end{abstract}

\maketitle

\section{Introduction}\label{sec:introduction}

Let $s_2(n)$ be the number of ones in the binary expansion of a nonnegative
integer $n$.  For $t\ge0$ and $d\in\Z$, B\'esineau \cite{Besineau1972}
proved the existence of the densities
\begin{equation}\label{eq:mu-definition}
 \mu_t(d)=\dens\{n\ge0:s_2(n+t)-s_2(n)=d\}.
\end{equation}
Here, for a set $A\subseteq\Z_{\ge0}$, $\dens(A)$ denotes its natural density
whenever it exists.  The measures $\mu_t$ are probability measures on $\Z$.
Cusick's conjecture concerns the half-line mass
\begin{equation}\label{eq:c-t-def}
 c_t:=\dens\{n\ge0:s_2(n+t)\ge s_2(n)\}=\mu_t(\Z_{\ge0}).
\end{equation}
The equality follows because $s_2(n+t)-s_2(n)\le s_2(t)$, so the event on the
left is a finite union of the level sets in \eqref{eq:mu-definition}.  The
conjecture, formulated in this form in the work of Drmota, Kauers and
Spiegelhofer \cite{DrmotaKauersSpiegelhofer2016} and attributed there to
T.~W. Cusick, asserts that
\[
 c_t>\frac12\qquad (t\ge1).
\]

The question is elementary to state but delicate because it is a uniform
problem about carry propagation.  Adding a fixed integer $t$ can create carry
chains through arbitrarily long blocks of binary digits, and these chains have
a nonlocal effect on $s_2$.  Cusick's conjecture says that, for every binary
pattern of $t$, this process has a strict sign bias: addition by $t$ weakly
increases the number of ones more often than it decreases it.  The word
``every'' is essential.  Typical or many-block patterns can be treated by
asymptotic methods, but sparse block patterns with very long blocks do not
fall under a normal-approximation regime.

The conjecture is connected with several equivalent or related problems.
Legendre's formula gives
\[
 s_2(n+t)-s_2(n)=s_2(t)-\nu_2\binom{n+t}{t},
\]
so the same density may be phrased in terms of the two-adic valuation of
binomial coefficients.  Drmota, Kauers and Spiegelhofer
\cite{DrmotaKauersSpiegelhofer2016} gave further equivalent formulations in
terms of carries, hyperbinary expansions, and zeros of
$X(X+1)\cdots(X+t-1)$ modulo $2^{t+1}$.  The problem is also related to the
Tu--Deng conjecture on binary strings, which arose in the construction of
Boolean functions with optimal algebraic immunity
\cite{TuDeng2011,CusickLiStanica2011,SpiegelhoferWallner2019}.

The threshold $1/2$ is critical.  For every $\varepsilon>0$, Drmota, Kauers
and Spiegelhofer proved that
\[
 \frac12<c_t<\frac12+\varepsilon
\]
for a set of parameters $t$ of asymptotic density one
\cite{DrmotaKauersSpiegelhofer2016}.  Thus no positive uniform gap above
$1/2$ can hold.  Subsequent work established the conjecture, or nearly so,
under hypotheses forcing many alternating blocks in the binary expansion of
$t$ \cite{Spiegelhofer2022,SpiegelhoferWallner2023}, and Sobolewski and
Spiegelhofer obtained a structural decomposition of the correlation measure
\cite{SobolewskiSpiegelhofer2026}.  These results leave open the genuinely
finite obstruction: the strict sign must be recovered for each individual
binary word, including words with very few blocks.

A finite form of the obstruction was already visible in
\cite{DrmotaKauersSpiegelhofer2016}.  Their simplified array, denoted below by
$P_t$, gives only the weak pair of inequalities
\[
 \widetilde c_t\le\frac12\le c_t,
 \qquad
 \widetilde c_t:=\dens\{n\ge0:s_2(n+t)>s_2(n)\},
\]
unless one can prove and transfer a strict finite bias. Tar\l{}owski
\cite{Tarlowski2026} subsequently introduced a stopped-simple-random-walk
representation of these finite measures, together with a recursively grown
binary-tree model, and established several basic identities and stabilization
properties for this framework. In particular, the persistence of the
one-sided asymmetry required for Cusick's conjecture was formulated as
Conjecture~24 of \cite{Tarlowski2026}. We recall, in our notation, the parts
of this framework that are needed below, with explicit references to the
corresponding results of \cite{Tarlowski2026}. The main new ingredient of the
present paper is a direct combinatorial identification of the recursively
defined trees with first-exit trees associated with principal subsequence
ideals. Under this identification, the required asymmetry becomes a finite
deletion problem. A marked-deletion argument then yields a quantitative
strict bias, which ultimately leads to the proof of Cusick's conjecture.

We prove the conjecture with an explicit nonuniform margin.

\begin{thm}\label{thm:main}
For every positive integer $t$,
\begin{equation}\label{eq:main-bound}
 c_t\ge \frac12+2^{-2s_2(t)-1}.
\end{equation}
In particular, $c_t>1/2$ for every $t\ge1$.
\end{thm}

The proof has four steps.  First, an exact deconvolution separates the
infinite trailing-one law of adding $1$ from a finitely supported probability
measure:
\[
 \mu_t=\mu_1*P_t.
\]
Second, if $t_w=(1w1)_2$ for a finite binary word $w$, then
\[
 P_{t_w}=\Law(S_{\tau_{\rev w}}),
\]
where $\tau_a$ is the first time that a random binary word leaves the
principal subsequence ideal $\Sub(a)$.  Thus the recursive tree for $P_{t_w}$
is exactly a first-exit tree for $\Sub(\rev w)$.  The reversal appears because
the integer recursion appends letters to the right of $w$, whereas the
subsequence-tree recursion prepends letters to the generator.

Third, we prove a median theorem for every finite deletion-closed binary
language.  Its proof is an exact level-by-level identity: the half-line
median defect is a weighted sum of integer defects comparing balanced words
with words having one more zero than one.  A marked-deletion double count
makes each defect nonnegative.  In a principal ideal, one terminal marked
insertion slot is missing, which forces a positive defect of size at least one
at a specified level.  This gives the quantitative first-exit bias.

Finally, a finite constant-tail expansion transfers the stopped-tree bias
back to the original digit-sum density.  Total-variation convergence alone
would not preserve a strict gap above $1/2$.  Instead, we prove that at the
explicit depth $s_2(t)$ the finite remainder and the omitted tail are already
supported below zero.  Hence the nonnegative half-line probability has
stabilized at a finite tree, and the strict first-exit defect yields
\eqref{eq:main-bound}.

The paper is organized as follows.  Section~\ref{sec:words-trees} fixes the
word, stopping-time, and binary-tree conventions.  Section~\ref{sec:correlation-trees}
recalls the correlation recurrences, defines the finite measures $P_t$, and
identifies their trees with principal subsequence ideals.  Sections
\ref{sec:median} and \ref{sec:strict} prove the weak and strict median
theorems for deletion-closed languages.  Section~\ref{sec:cusick} combines
the finite constant-tail expansion with the support cutoff and proves
Theorem~\ref{thm:main}.

\section{Words, stopped walks, and binary trees}\label{sec:words-trees}

\subsection{Deletion-closed languages}
Let $\bits$ be the set of all finite binary words, including the empty word
$\eps$.  For $u\in\bits$, write $|u|$ for its length, $|u|_0$ and $|u|_1$ for the numbers of zeros and ones, and
\begin{equation*}
 b(u)=|u|_1-|u|_0
\end{equation*}
for its balance.  Concatenation is written without a symbol.

\begin{defn}\label{def:deletion-closed}
A finite nonempty set $D\subseteq\bits$ is called a \emph{deletion-closed language}, or a \emph{subsequence order ideal}, if every word obtained from an
element of $D$ by deleting letters also belongs to $D$.
\end{defn}

Throughout this paper, the word ``ideal'' means a down-set in the subsequence order; no directedness assumption is intended. Every deletion-closed language
contains $\eps$ and is prefix-closed. For a fixed word $a$, the principal order ideal is
\begin{equation*}
 \Sub(a)=\{u\in\bits:u\preccurlyeq a\}.
\end{equation*}

Throughout the probabilistic parts of the paper, let
$\xi_1,\xi_2,\ldots$ be independent binary random variables with
\[
 \mathbb P(\xi_j=0)=\mathbb P(\xi_j=1)=\frac12,
\]
and let
\[
 S_0=0,\qquad S_n=\sum_{j=1}^{n}(2\xi_j-1)
\]
be the simple symmetric random walk obtained by reading a $1$ as a step $+1$
and a $0$ as a step $-1$.  For a word $a$, set
\begin{equation}\label{eq:tau-a-intro}
 \tau_a=\min\{n\ge1:\xi_1\ldots\xi_n\not\preccurlyeq a\}.
\end{equation}
Equivalently, the sampled word is read until its current prefix first leaves
$\Sub(a)$.

For a finite deletion-closed language $D$, define
\begin{equation*}
 \tau_D=\min\{n\ge1:\xi_1\ldots\xi_n\notin D\}.
\end{equation*}
If $L=\max\{|u|:u\in D\}$, then $\tau_D\le L+1$. Thus every expectation below
is a finite sum. For $D=\Sub(a)$, the notation $\tau_D$ agrees with
$\tau_a$ in \eqref{eq:tau-a-intro}.

\subsection{Full binary trees and stopped laws}
A finite full rooted binary tree is either a leaf, denoted by $\bullet$, or an
ordered pair $[T_0,T_1]$ of full binary trees.  The edges to $T_0$ and $T_1$ are labelled $0$ and $1$, respectively.  Define a probability measure $E(T)$ on $\Z$ recursively by
\begin{equation}\label{eq:tree-law}
 E(\bullet)=\delta_0,
 \qquad
E([T_0,T_1])=\frac12\sigma_{-1}E(T_0)+\frac12\sigma_1E(T_1),
\end{equation}
where the translation operator is
\begin{equation*}
 (\sigma_j\nu)(d)=\nu(d-j).
\end{equation*}
Thus $\sigma_j\nu$ is the law of $X+j$ when $\nu=\Law(X)$.

If $D$ is finite and prefix-closed, let $\mathcal T(D)$ be the full binary tree
whose internal nodes, identified with root-to-node words, are exactly the words in $D$.  A path stops at the first prefix outside $D$, and hence
\begin{equation}\label{eq:tree-stopped-law}
 E(\mathcal T(D))=\Law(S_{\tau_D}).
\end{equation}
For a word $a$, write
\begin{equation*}
 \mathcal S(a)=\mathcal T(\Sub(a)).
\end{equation*}
If $T=[T_0,T_1]$, we write $T^{(0)}=T_0$ and $T^{(1)}=T_1$.

\section{Correlation measures and principal-subsequence trees}
\label{sec:correlation-trees}

\subsection{Correlation recurrences and deconvolution}
Put
\[
 d_t(n)=s_2(n+t)-s_2(n).
\]
The identities $s_2(2n)=s_2(n)$ and $s_2(2n+1)=s_2(n)+1$ give, for
$n,t\ge0$,
\begin{align*}
 d_{2t}(2n)&=d_t(n),
 &d_{2t}(2n+1)&=d_t(n),\\
 d_{2t+1}(2n)&=d_t(n)+1,
 &d_{2t+1}(2n+1)&=d_{t+1}(n)-1.
\end{align*}
Splitting the nonnegative integers according to parity therefore yields the standard binary recurrences
\begin{align}
\mu_{2t}&=\mu_t,\label{eq:mu-even}\\
\mu_{2t+1}&=\frac12\sigma_{-1}\mu_{t+1}+\frac12\sigma_1\mu_t.\notag
\end{align}
where $\mu_t$ is the probability measure defined in \eqref{eq:mu-definition}; here $t\ge0$, and $\mu_0=\delta_0$.  Finally, if the binary expansion of $n$
has exactly $j$ trailing ones, then $d_1(n)=1-j$; this event has density
$2^{-j-1}$. Hence
\begin{equation}\label{eq:mu1}
\mu_1=\sum_{j=0}^{\infty}2^{-j-1}\delta_{1-j}.
\end{equation}

For $t\ge1$, define a finitely supported probability measure $P_t$ by
\begin{align}
P_1&=\delta_0,\label{eq:P1}\\
P_{2t}&=P_t,\label{eq:P-even}\\
P_{2t+1}&=\frac12\sigma_{-1}P_{t+1}+\frac12\sigma_1P_t\qquad(t\ge1).\label{eq:P-odd}
\end{align}
These are the finite deconvolution measures appearing in
\cite{DrmotaKauersSpiegelhofer2016,Tarlowski2026}.  In the notation of
Drmota--Kauers--Spiegelhofer, one has
\begin{equation}\label{eq:P-varphi}
P_t(d)=\varphi(d,t),
\end{equation}
so their half-line quantity $p_t$ is $P_t(\Z_{\ge0})$.
The two identities in the following lemma appear, in the corresponding
$L/R$ notation, as Corollary~7 and equation~(2.12), respectively, in
Tar\l{}owski~\cite{Tarlowski2026}.  We include a short proof to keep the
normalizations and shift conventions explicit.

\begin{lem}\label{lem:convolution}
For every $t\ge1$,
\begin{equation}\label{eq:convolution}
 \mu_t=\mu_1*P_t.
\end{equation}
Furthermore, for every $d\in\Z$,
\begin{equation}\label{eq:deconvolution}
 P_t(d)=2\mu_t(d+1)-\mu_t(d+2).
\end{equation}
\end{lem}

\begin{proof}
We first prove \eqref{eq:convolution} by strong induction on $t$.  The case
$t=1$ follows from $P_1=\delta_0$. The even step follows from
\eqref{eq:mu-even} and \eqref{eq:P-even}.  If $t=2s+1\ge3$, then convolution commutes with translation, and the induction hypothesis gives
\begin{align*}
 \mu_{2s+1}
 &=\frac12\sigma_{-1}(\mu_1*P_{s+1})
+\frac12\sigma_1(\mu_1*P_s)\\
&=\mu_1*\left(\frac12\sigma_{-1}P_{s+1}+\frac12\sigma_1P_s\right)=\mu_1*P_{2s+1}.
\end{align*}

From \eqref{eq:mu1}, one checks directly that
\begin{equation}\label{eq:mu1-inverse-kernel}
2\mu_1(r+1)-\mu_1(r+2)=\one_{\{r=0\}}\qquad(r\in\Z).
\end{equation}
Applying \eqref{eq:mu1-inverse-kernel} term by term to
$\mu_t=\mu_1*P_t$ gives
\begin{align*}
 2\mu_t(d+1)-\mu_t(d+2)
 &=\sum_{x\in\Z}P_t(x)
\bigl(2\mu_1(d+1-x)-\mu_1(d+2-x)\bigr)\\
&=\sum_{x\in\Z}P_t(x)\one_{\{d=x\}}=P_t(d),
\end{align*}
which is \eqref{eq:deconvolution}.
\end{proof}

\subsection{The recursive trees}
The following construction is the $0/1$-indexed form of the growing-tree
construction in Section~3.2 of Tar\l{}owski~\cite{Tarlowski2026}.  Under the
dictionary $0\leftrightarrow L$ and $1\leftrightarrow R$, his recursion
(3.1) gives the recursion below on odd indices; the convention $T_{2t}=T_t$
merely removes powers of two.  His Theorem~11 gives the corresponding stopped
law $E(T_t)=P_t$.  We restate these facts in our notation because the new step
used later is the subsequent identification of these recursively grown trees
with principal-subsequence first-exit trees.

Define full binary trees $T_t$ for $t\ge1$ by
\begin{equation}\label{eq:T-recursion}
T_1=\bullet,\qquad T_{2t}=T_t,\qquad
T_{2t+1}=[T_{t+1},T_t]\quad(t\ge1).
\end{equation}
Equations \eqref{eq:tree-law} and \eqref{eq:P1}--\eqref{eq:P-odd} imply
\begin{equation}\label{eq:E-P}E(T_t)=P_t.
\end{equation}

The next elementary recursion is the bridge from these number-theoretic trees to subsequence order ideals.

\begin{lem}\label{lem:subsequence-tree-recursion}
For every finite binary word $a$,
\begin{align}
\mathcal S(0a)&=[\mathcal S(a),\mathcal S(a)^{(1)}],\label{eq:S0a}\\
 \mathcal S(1a)&=[\mathcal S(a)^{(0)},\mathcal S(a)].\label{eq:S1a}
\end{align}
\end{lem}

\begin{proof}
For a word $c$ and a first letter $b\in\{0,1\}$, the internal-node language in
the subtree reached after reading $b$ is
\begin{equation*}
 D_b(c)=\{u\in\bits:bu\preccurlyeq c\}.
\end{equation*}
Consider $c=0a$. If $u\preccurlyeq a$, then $0u\preccurlyeq0a$ by matching
the first zero of $0u$ with the initial zero of $0a$. Conversely, in any embedding of $0u$ into $0a$, all letters realizing $u$ occur after the position realizing the initial zero of $0u$; hence they all lie among the letters of $a$, and therefore $u\preccurlyeq a$. Hence
\[
 D_0(0a)=\Sub(a).
\]
On the other hand, the initial zero of $0a$ cannot match a first letter $1$, so
\[
 D_1(0a)=D_1(a),
\]
which is exactly the internal-node language of the right subtree of
$\mathcal S(a)$.  This proves \eqref{eq:S0a}.  Interchanging zero and one proves
\eqref{eq:S1a}.
\end{proof}

For $w\in\bits$, recall $t_w={(1w1)}_2$.  Then
\begin{equation}\label{eq:t-word-recursion}
 t_{\eps}=3,
 \qquad
 t_{w0}=2t_w-1,
 \qquad
 t_{w1}=2t_w+1.
\end{equation}

\begin{prop}
\label{prop:subsequence-representation}
For every $w\in\bits$,
\begin{equation}\label{eq:tree-identification}
 T_{t_w}=\mathcal S(\rev w).
\end{equation}
Consequently,
\begin{equation}\label{eq:P-stopped-law}
P_{t_w}=\Law(S_{\tau_{\rev w}}).
\end{equation}
\end{prop}

\begin{proof}
For $w=\eps$, the recursion \eqref{eq:T-recursion} gives
\[
 T_3=[T_2,T_1]=[\bullet,\bullet]=\mathcal S(\eps).
\]
Assume \eqref{eq:tree-identification} holds for a word $w$, and write the odd integer $t_w=2u+1$.  Then
\begin{equation*}
 T_{t_w}=[T_{u+1},T_u].
\end{equation*}
Using \eqref{eq:T-recursion} and \eqref{eq:t-word-recursion}, we obtain
\begin{align}
 T_{t_{w0}}
 &=T_{2t_w-1}=T_{4u+1}
   =[T_{2u+1},T_{2u}]
   =[T_{t_w},T_u]
   =[T_{t_w},T_{t_w}^{(1)}],
 \label{eq:T-w0}\\
 T_{t_{w1}}
 &=T_{2t_w+1}=T_{4u+3}
   =[T_{2u+2},T_{2u+1}]
   =[T_{u+1},T_{t_w}]
   =[T_{t_w}^{(0)},T_{t_w}].
 \label{eq:T-w1}
\end{align}
Since
\[
 \rev{(w0)}=0\rev w,
 \qquad
 \rev{(w1)}=1\rev w,
\]
the induction step follows from Lemma~\ref{lem:subsequence-tree-recursion}.
Finally, \eqref{eq:P-stopped-law} follows from
\eqref{eq:E-P}, \eqref{eq:tree-stopped-law}, and
\eqref{eq:tree-identification}.
\end{proof}

\begin{exa}\label{exa:w01}
Take $w=01$.  Then $t_w={(1011)}_2=11$ and $\rev w=10$. The
internal-node language is
\[
 \Sub(10)=\{\eps,0,1,10\}.
\]
The stopping words are
\[
 00,\quad 01,\quad 11,\quad 100,\quad 101,
\]
with probabilities $1/4,1/4,1/4,1/8,1/8$, respectively.  Their balances are
$-2,0,2,-1,1$.  Therefore
\begin{equation*}
P_{11}=\frac14\delta_{-2}+\frac18\delta_{-1}
+\frac14\delta_0+\frac18\delta_1+\frac14\delta_2.
\end{equation*}
This example makes both the reversal in
\eqref{eq:tree-identification} and the path-probability convention explicit.
\end{exa}

\begin{rem}
\label{rem:dictionary}
Tar\l{}owski's operations are $L(t)=2t-1$ and $R(t)=2t+1$. By
\eqref{eq:t-word-recursion}, appending a digit $0$ to $w$ applies $L$ to $t_w$, while appending a digit $1$ applies $R$. The reversal in \eqref{eq:tree-identification} occurs because these operations append a digit to the integer word, whereas Lemma~\ref{lem:subsequence-tree-recursion}
prepends a letter to the principal generator.
\end{rem}

\section{An exact median formula and a deletion double count}
\label{sec:median}

Let $D$ be a finite deletion-closed language.  For $r\in\Z$, define the
weighted balance-level sum
\begin{equation*}
A_r(D)=\sum_{\substack{u\in D\\b(u)=r}}2^{-|u|}.
\end{equation*}

\begin{lem}\label{lem:stopped-tree-identity}
For every function $f:\Z\to\mathbb R$,
\begin{equation}\label{eq:stopped-tree-identity}
 \mathbb E f(S_{\tau_D})
 =f(0)+\sum_{u\in D}2^{-|u|}\Delta f(b(u)),
\end{equation}
where
\begin{equation*}
 \Delta f(x)=\frac{f(x-1)+f(x+1)}2-f(x).
\end{equation*}
\end{lem}

\begin{proof}
Let $L=\max\{|u|:u\in D\}$.  Since $\tau_D\le L+1$, telescoping gives
\begin{equation}\label{eq:telescoping-stopped}
 f(S_{\tau_D})-f(0)
 =\sum_{n=0}^{L}\one_{\{\tau_D>n\}}
   \bigl(f(S_{n+1})-f(S_n)\bigr).
\end{equation}
Conditioning on the first $n$ bits, and using the fairness of $\xi_{n+1}$,
the conditional expectation of the increment in
\eqref{eq:telescoping-stopped} is $\Delta f(S_n)$.  Because $D$ is
prefix-closed,
\[
 \{\tau_D>n\}=\{\xi_1\ldots\xi_n\in D\}.
\]
For each $u\in D$ of length $n$, the cylinder event
$\{\xi_1\ldots\xi_n=u\}$ has probability $2^{-|u|}$ and satisfies
$S_n=b(u)$.  Taking expectations in \eqref{eq:telescoping-stopped} and summing
over all cylinders gives \eqref{eq:stopped-tree-identity}.
\end{proof}

For $m\ge1$, introduce the three finite level sets
\begin{align*}
 \mathcal B_m(D)
 &=\{u\in D:|u|_0=|u|_1=m\},\\
 \mathcal C_m^-(D)
 &=\{v\in D:|v|_0=m,\ |v|_1=m-1\},\\
 \mathcal C_m^+(D)
 &=\{v\in D:|v|_0=m-1,\ |v|_1=m\}.
\end{align*}
Define the integer level defects
\begin{equation*}
 \Delta_m^-(D)=2|\mathcal C_m^-(D)|-|\mathcal B_m(D)|,
 \qquad
 \Delta_m^+(D)=2|\mathcal C_m^+(D)|-|\mathcal B_m(D)|.
\end{equation*}

\begin{prop}\label{prop:exact-median-defect}
For every finite deletion-closed language $D$,
\begin{align}
 \mathbb P(S_{\tau_D}\ge0)
&=\frac12+\sum_{m\ge1}2^{-2m-1}\Delta_m^-(D),
 \label{eq:exact-plus}\\
 \mathbb P(S_{\tau_D}\le0)
&=\frac12+\sum_{m\ge1}2^{-2m-1}\Delta_m^+(D).
 \label{eq:exact-minus}
\end{align}
Both sums are finite.
\end{prop}

\begin{proof}
Let $f_+(x)=\one_{\{x\ge0\}}$.  Its discrete Laplacian is zero except at
$-1$ and $0$, where
\[
 \Delta f_+(-1)=\frac12,
 \qquad
 \Delta f_+(0)=-\frac12.
\]
Since $f_+(0)=1$, Lemma~\ref{lem:stopped-tree-identity} gives
\begin{equation}\label{eq:prob-plus-A}
 \mathbb P(S_{\tau_D}\ge0)
 =1+\frac12A_{-1}(D)-\frac12A_0(D)
=\frac12+\frac12\bigl(A_{-1}(D)-(A_0(D)-1)\bigr).
\end{equation}
Now
\begin{align}
 A_0(D)-1
&=\sum_{m\ge1}2^{-2m}|\mathcal B_m(D)|,
 \label{eq:A0-levels}\\
 A_{-1}(D)
 &=\sum_{m\ge1}2^{-(2m-1)}|\mathcal C_m^-(D)|
=\sum_{m\ge1}2^{-2m}\,2|\mathcal C_m^-(D)|.
 \label{eq:Aminus1-levels}
\end{align}
Substituting \eqref{eq:A0-levels} and \eqref{eq:Aminus1-levels} into
\eqref{eq:prob-plus-A} yields \eqref{eq:exact-plus}.

For $f_-(x)=\one_{\{x\le0\}}$, its discrete Laplacian is zero except at
$0$ and $1$, where
\[
 \Delta f_-(0)=-\frac12,
 \qquad
 \Delta f_-(1)=\frac12.
\]
Lemma~\ref{lem:stopped-tree-identity} therefore gives
\[
 \mathbb P(S_{\tau_D}\le0)
 =1+\frac12A_1(D)-\frac12A_0(D).
\]
Moreover,
\begin{equation}\label{eq:Aplus1-levels}
 A_1(D)
 =\sum_{m\ge1}2^{-(2m-1)}|\mathcal C_m^+(D)|
=\sum_{m\ge1}2^{-2m}\,2|\mathcal C_m^+(D)|.
\end{equation}
Combining this identity with \eqref{eq:A0-levels} yields
\eqref{eq:exact-minus}.
\end{proof}

The next lemma proves that every summand in
\eqref{eq:exact-plus} and \eqref{eq:exact-minus} is nonnegative.  We spell out
the marked objects because two different insertion slots may produce the same
unmarked word.

\begin{lem}
\label{lem:deletion-inequality}
For every finite deletion-closed language $D$ and every $m\ge1$,
\begin{equation*}
 \Delta_m^-(D)\ge0,
 \qquad
 \Delta_m^+(D)\ge0.
\end{equation*}
Equivalently,
\[
 |\mathcal B_m(D)|\le2|\mathcal C_m^-(D)|,
 \qquad
 |\mathcal B_m(D)|\le2|\mathcal C_m^+(D)|.
\]
\end{lem}

\begin{proof}
We prove the first inequality.  Let
\begin{equation*}
 \mathcal M_m^-(D)
 =\{(u,i):u\in\mathcal B_m(D),\ 1\le i\le2m,\ u_i=1\}.
\end{equation*}
Thus the second coordinate marks one of the $m$ ones of $u$, and
\begin{equation}\label{eq:marked-source-count}
 |\mathcal M_m^-(D)|=m|\mathcal B_m(D)|.
\end{equation}

If $v=v_1\cdots v_{2m-1}$, it has exactly $2m$ insertion slots, indexed by
$j=0,1,\ldots,2m-1$: slot $j$ inserts a new letter after the first $j$
letters.  Let $I_j^1(v)$ be the word obtained by inserting a $1$ in slot $j$, and define the set of admissible marked slots
\begin{equation}\label{eq:admissible-slots}
 \mathcal I_m^-(D)
 =\{(v,j):v\in\mathcal C_m^-(D),\ 0\le j\le2m-1,
\ I_j^1(v)\in D\}.
\end{equation}
Deleting the marked one from $(u,i)\in\mathcal M_m^-(D)$ produces a word
$v\in\mathcal C_m^-(D)$ because $D$ is deletion-closed, and records the slot $j=i-1$. Conversely, inserting a marked one in an admissible slot
$(v,j)$ recovers the unique source pair $(I_j^1(v),j+1)$.  Hence marked deletion is a bijection
\begin{equation}\label{eq:marked-bijection}
 \mathcal M_m^-(D)\longleftrightarrow\mathcal I_m^-(D).
\end{equation}
The same unmarked word can occur for two adjacent slots when $v$ contains
ones, but the inserted letter is marked, so the corresponding source pairs are distinct. From \eqref{eq:marked-source-count} and the trivial bound
\[
 |\mathcal I_m^-(D)|\le2m|\mathcal C_m^-(D)|
\]
we obtain
\[
 m|\mathcal B_m(D)|\le2m|\mathcal C_m^-(D)|.
\]
Dividing by $m$ gives $\Delta_m^-(D)\ge0$.  Deleting a marked zero and
inserting zero into a word of $\mathcal C_m^+(D)$ proves the second inequality.
\end{proof}

\begin{thm}
\label{thm:ideal-median}
For every finite deletion-closed language $D$,
\begin{equation}\label{eq:ideal-median}
 \mathbb P(S_{\tau_D}\ge0)\ge\frac12,
 \qquad
 \mathbb P(S_{\tau_D}\le0)\ge\frac12.
\end{equation}
\end{thm}

\begin{proof}
Every level defect in \eqref{eq:exact-plus} and \eqref{eq:exact-minus} is
nonnegative by Lemma~\ref{lem:deletion-inequality}.  The conclusion follows
immediately from Proposition~\ref{prop:exact-median-defect}.
\end{proof}

\begin{rem}\label{rem:not-unmarked-injection}
The proof of Lemma~\ref{lem:deletion-inequality} is best viewed as a double
count, or equivalently as a bijection between marked deletions and admissible marked insertion slots.  Describing it merely as an injection between
unmarked word sets would hide the multiplicity caused by consecutive equal letters.
\end{rem}

\section{Strict asymmetry for principal subsequence ideals}
\label{sec:strict}

The principal structure $D=\Sub(a)$ guarantees that at least one of the
insertion slots counted above is not admissible.

\begin{lem}
\label{lem:unfilled-slot}
Let $a$ be a nonempty binary word ending in $0$.  Suppose that $a$ has $Z$ zeros and $O$ ones, and put
\begin{equation*}
 m=\min(Z,O+1).
\end{equation*}
There exists a word $q\in\Sub(a)$ such that
\begin{equation*}
 |q|_0=m,
 \qquad
 |q|_1=m-1,
 \qquad
 q1\notin\Sub(a).
\end{equation*}
\end{lem}

\begin{proof}
There are two cases.

Suppose first that $Z\ge O+1$.  Retain, in their original order, all $O$ ones
of $a$ and any $O+1$ zeros.  The resulting subsequence $q$ has
$|q|_0=O+1=m$ and $|q|_1=O=m-1$.  The word $q1$ contains $O+1$ ones, whereas
$a$ contains only $O$, so $q1\not\preccurlyeq a$.

Now suppose that $Z\le O$.  Retain all $Z$ zeros of $a$ and any $Z-1$ of its ones.  Then $q\preccurlyeq a$, $|q|_0=Z=m$, and $|q|_1=Z-1=m-1$.  We claim
that an embedding of $q$ into $a$ must map the $i$th zero of $q$ to the $i$th zero of $a$ for every $i$. Indeed, the embedding is order preserving and $q$ and $a$ contain the same number $Z$ of zeros; hence the induced increasing map between their zero positions is a bijection, and therefore is the order-preserving bijection.  In particular, the last zero of $q$ must map to the last zero of $a$. Since $a$ ends in that zero, no position remains in
which to place the appended $1$. Equivalently, the retained ones of $q$ all occur before that last zero of $a$, while the final zero of $q$ is forced to occupy the final position of $a$; hence the appended $1$ in $q1$ would have to be embedded after position $|a|$, which is impossible. Thus $q1\not\preccurlyeq a$.
\end{proof}

\begin{exa}\label{exa:missing-slot}
Take $a=11010$. Then $Z=2$, $O=3$, $m=2$, and the second case of
Lemma~\ref{lem:unfilled-slot} permits $q=100$.  The four possible marked
insertions of a $1$ into $q$ are
\[
 \widehat{1}100,
 \qquad
 1\widehat{1}00,
 \qquad
 10\widehat{1}0,
 \qquad
 100\widehat{1}.
\]
The first two marked objects have the same underlying word $1100$, but they are different source pairs. The terminal insertion has underlying word $1001$, which is not a subsequence of $11010$.  Thus at least one of the
$2m=4$ marked slots over the target $q$ is absent.
\end{exa}

\begin{thm}
\label{thm:strict-principal}
Let $a$ be a nonempty binary word.
\begin{enumerate}[label=\textup{(\roman*)},leftmargin=2.7em]
\item If $a$ ends in $0$, has $Z$ zeros and $O$ ones, and $m=\min(Z,O+1)$, then
\begin{equation*}
 \mathbb P(S_{\tau_a}\ge0)
 \ge \frac12+2^{-2m-1}.
\end{equation*}
\item If $a$ ends in $1$, has $Z$ zeros and $O$ ones, and
$m=\min(O,Z+1)$, then
\begin{equation*}
 \mathbb P(S_{\tau_a}\le0)
 \ge \frac12+2^{-2m-1}.
\end{equation*}
\end{enumerate}
\end{thm}

\begin{proof}
Assume first that $a$ ends in $0$, and set $D=\Sub(a)$. Let
$q\in\mathcal C_m^-(D)$ be supplied by Lemma~\ref{lem:unfilled-slot}.  In the admissible-slot set
\eqref{eq:admissible-slots}, the terminal slot $(q,2m-1)$ is absent because
$I_{2m-1}^1(q)=q1\notin D$. Therefore the marked bijection
\eqref{eq:marked-bijection} gives the strict count
\begin{equation*}
 m|\mathcal B_m(D)|
 =|\mathcal I_m^-(D)|
 \le2m|\mathcal C_m^-(D)|-1.
\end{equation*}
Rearranging,
\begin{equation}\label{eq:multiple-of-m}
 m\Delta_m^-(D)
 =m\bigl(2|\mathcal C_m^-(D)|-|\mathcal B_m(D)|\bigr)\ge1.
\end{equation}
The left-hand side of \eqref{eq:multiple-of-m} is a multiple of the positive integer $m$. Hence it is at least $m$, and
\begin{equation}\label{eq:strict-level-defect}
 \Delta_m^-(D)\ge1.
\end{equation}
All other level defects are nonnegative by
Lemma~\ref{lem:deletion-inequality}. Substituting
\eqref{eq:strict-level-defect} into the exact formula
\eqref{eq:exact-plus} yields
\[
 \mathbb P(S_{\tau_a}\ge0)
 \ge\frac12+2^{-2m-1}.
\]

For part (ii), complement every bit of $a$ and of every sampled path.  The
principal ideal is carried bijectively to the principal ideal of the
complemented generator, zeros and ones are interchanged, and the stopped position changes sign. Part (ii) therefore follows from part (i).
\end{proof}

\begin{thm}\label{thm:subsequence-stopping}
For every finite binary word $a$,
\[
 \mathbb P(S_{\tau_a}\ge0)\ge\frac12,
 \qquad
 \mathbb P(S_{\tau_a}\le0)\ge\frac12.
\]
If $a$ is nonempty and ends in $0$, has $Z$ zeros and $O$ ones, and
$m=\min(Z,O+1)$, then
\[
 \mathbb P(S_{\tau_a}\ge0)
 \ge \frac12+2^{-2m-1}.
\]
If $a$ ends in $1$, has $Z$ zeros and $O$ ones, and
$m=\min(O,Z+1)$, then
\[
 \mathbb P(S_{\tau_a}\le0)
 \ge \frac12+2^{-2m-1}.
\]
\end{thm}

\begin{proof}
The two weak median inequalities follow from Theorem~\ref{thm:ideal-median}
applied to $D=\Sub(a)$.  The two quantitative strict inequalities are exactly
Theorem~\ref{thm:strict-principal}.
\end{proof}

Combining Theorems~\ref{thm:ideal-median} and \ref{thm:strict-principal} with
Proposition~\ref{prop:subsequence-representation} gives the tree statement
needed later.

\begin{cor}\label{cor:tree-asymmetry}
Let $w\in\bits$. Then
\begin{equation*}
 P_{t_w}(\Z_{\ge0})\ge\frac12,
 \qquad
 P_{t_w}(\Z_{\le0})\ge\frac12.
\end{equation*}
If $w$ is nonempty and begins with $0$, then
\begin{equation}\label{eq:quantitative-tree-plus}
 P_{t_w}(\Z_{\ge0})
 \ge\frac12+
2^{-2\min(|w|_0,|w|_1+1)-1}.
\end{equation}
If $w$ is nonempty and begins with $1$, then
\begin{equation}\label{eq:quantitative-tree-minus}
 P_{t_w}(\Z_{\le0})
\ge\frac12+2^{-2\min(|w|_1,|w|_0+1)-1}.
\end{equation}
\end{cor}

\begin{proof}
By Proposition~\ref{prop:subsequence-representation}, $P_{t_w}$ is the stopped law associated with $\Sub(\rev w)$. The first letter of $w$ is the
last letter of $\rev w$, while reversal preserves the numbers of zeros and
ones. The weak assertions follow from Theorem~\ref{thm:ideal-median}, and
\eqref{eq:quantitative-tree-plus}--\eqref{eq:quantitative-tree-minus} follow
from the two parts of Theorem~\ref{thm:strict-principal}.
\end{proof}

\begin{rem}\label{rem:dks-halfline}
By \eqref{eq:P-varphi}, the first inequality in
Corollary~\ref{cor:tree-asymmetry}, together with $P_1=\delta_0$ and the even recurrence \eqref{eq:P-even}, proves
\[
p_t=P_t(\Z_{\ge0})\ge\frac12\qquad(t\ge1),
\]
the conjectural half-line inequality isolated in
\cite{DrmotaKauersSpiegelhofer2016}.  Their Lemma~3.1 then yields the weak pair
$\widetilde c_t\le1/2\le c_t$.  The strict step requires the principal-ideal
asymmetry developed here, together with the finite-stabilization mechanism that
already appears in the proof of Lemma~25 of Tar\l{}owski~\cite{Tarlowski2026}.
\end{rem}

\begin{rem}\label{rem:tarlowski-conjecture}
Under the dictionary $0\leftrightarrow L$ and $1\leftrightarrow R$ from
Remark~\ref{rem:dictionary}, Corollary~\ref{cor:tree-asymmetry} proves the
strict persistence assertion formulated as Conjecture~24 in
\cite{Tarlowski2026}.  The estimates
\eqref{eq:quantitative-tree-plus}--\eqref{eq:quantitative-tree-minus} strengthen
that conjecture quantitatively for each word. The mechanism is the
principal-subsequence representation followed by the absent insertion slot.
\end{rem}

\section{Finite stabilization and the proof of Cusick's conjecture}
\label{sec:cusick}

In this section $0^n$ denotes the word consisting of $n$ zeros. The exact finite expansions and support cutoffs used below are proved directly from the recursions above, so the final transfer from stopped trees to $c_t$ is logically self-contained.

For readability, write
\begin{equation*}
 P_w:=P_{t_w},
 \qquad
 T_w:=T_{t_w}=\mathcal S(\rev w).
\end{equation*}
Let
\begin{equation*}
 P_w^{(0)}=E(T_w^{(0)})
\end{equation*}
be the law attached to the left subtree of $T_w$.  This is the subtree law
itself; it does not include the initial step $-1$ taken to enter the left
subtree. For two probability
measures $\nu$ and $\eta$ on $\Z$, we use the standard total-variation distance
\begin{equation}\label{eq:tv-definition}
 d_{\mathrm{TV}}(\nu,\eta)
 =\frac12\sum_{d\in\Z}|\nu(d)-\eta(d)|.
\end{equation}

\subsection{The finite comb identity}
The constant-tail limit is part of the tree dynamics developed in
\cite{Tarlowski2026}.  We give the exact finite expansion because both its
remainder and its support are used in the proof of
Theorem~\ref{thm:main}.

\begin{lem}\label{lem:finite-comb}
For every $w\in\bits$ and $n\ge0$,
\begin{equation}\label{eq:finite-comb}
 P_{w10^n}
 =2^{-n-1}\sigma_{-n-1}P_w^{(0)}
  +\sum_{j=0}^{n}2^{-j-1}\sigma_{1-j}P_w.
\end{equation}
Moreover,
\begin{equation}\label{eq:comb-tv-bound}
 d_{\mathrm{TV}}(P_{w10^n},\mu_{t_w})\le2^{-n-1}.
\end{equation}
In particular,
\begin{equation}\label{eq:comb-limit}
P_{w10^n}\longrightarrow\mu_{t_w}
 \quad\text{in total variation.}
\end{equation}
\end{lem}

\begin{proof}
Equation \eqref{eq:T-w1} gives
\begin{equation}\label{eq:first-comb-tree}
 T_{w1}=[T_w^{(0)},T_w].
\end{equation}
Equation \eqref{eq:T-w0} shows that appending a zero replaces a tree $T_x$ by
$[T_x,T_x^{(1)}]$.  Since the right subtree of $T_{w1}$ is $T_w$, induction on $n$ yields
\begin{equation}\label{eq:comb-tree-recursion}
 T_{w10^n}=[T_{w10^{n-1}},T_w]\qquad(n\ge1).
\end{equation}
Put $Q_n=P_{w10^n}$.  From \eqref{eq:tree-law},
\eqref{eq:first-comb-tree}, and \eqref{eq:comb-tree-recursion},
\begin{align}
Q_0&=\frac12\sigma_{-1}P_w^{(0)}+\frac12\sigma_1P_w,
 \notag\\
Q_n&=\frac12\sigma_{-1}Q_{n-1}+\frac12\sigma_1P_w
 \qquad(n\ge1).
 \label{eq:Q-recursion}
\end{align}
A direct induction in \eqref{eq:Q-recursion} gives
\eqref{eq:finite-comb}.

By \eqref{eq:mu1} and Lemma~\ref{lem:convolution},
\begin{equation}\label{eq:infinite-comb}
 \sum_{j=0}^{\infty}2^{-j-1}\sigma_{1-j}P_w
 =\mu_1*P_w
 =\mu_{t_w}.
\end{equation}
The first term on the right of \eqref{eq:finite-comb} has total mass
$2^{-n-1}$. The omitted tail in \eqref{eq:infinite-comb} has total mass
\[
 \sum_{j=n+1}^{\infty}2^{-j-1}=2^{-n-1}.
\]
Consequently,
\[
 \sum_{d\in\Z}|P_{w10^n}(d)-\mu_{t_w}(d)|\le2^{-n},
\]
and \eqref{eq:comb-tv-bound} follows from
\eqref{eq:tv-definition}.
\end{proof}

\begin{rem}
\label{rem:finite-comb-prior}
Let $v\in\{L,R\}^*$ be obtained from $w$ by replacing each $0$ with $L$
and each $1$ with $R$.  Under the dictionary of
Remark~\ref{rem:dictionary}, the word $w10^n$ corresponds to $vRL^n$.
Accordingly, formula \eqref{eq:finite-comb} is the $0/1$-notation version of
equation~(4.10) in Tar\l{}owski~\cite{Tarlowski2026}, and the limiting identity
\eqref{eq:comb-limit} is the corresponding case of his Theorem~21.  We have
included the finite induction to make every shift and coefficient auditable and
to record the explicit total-variation estimate \eqref{eq:comb-tv-bound}.
\end{rem}

\subsection{Support and reversal}

\begin{lem}\label{lem:support-bound}
For every $w\in\bits$,
\begin{equation*}
 \max\supp(P_w)=|w|_1+1,
 \qquad
\max\supp(P_w^{(0)})\le|w|_1+1.
\end{equation*}
\end{lem}

\begin{proof}
Put $a=\rev w$ and $r=|a|_1=|w|_1$. A stopping word $y$ for $\Sub(a)$ has its proper prefix $y^-$ in $\Sub(a)$. Hence $y^-$ contains at most $r$ ones,
so $y$ contains at most $r+1$ ones. Therefore
\[
 b(y)=|y|_1-|y|_0\le|y|_1\le r+1.
\]
The stopping word $1^{r+1}$ attains equality: $1^r\preccurlyeq a$, whereas
$1^{r+1}\not\preccurlyeq a$. This proves the first assertion.

If the left subtree is a leaf, then $P_w^{(0)}=\delta_0$ and the second
assertion is immediate.  Otherwise its internal-node language is
\[
 \{u\in\bits:0u\preccurlyeq a\}.
\]
If $y$ is a stopping word within that subtree, then $0y^-\preccurlyeq a$.
Thus $y^-$ contains at most $r$ ones and $y$ contains at most $r+1$ ones, which
gives the second assertion.
\end{proof}

\begin{lem}\label{lem:word-reversal}
For every finite binary word $a$,
\begin{equation}\label{eq:stopped-reversal}
\Law(S_{\tau_a})=\Law(S_{\tau_{\rev a}}).
\end{equation}
Consequently, for every $w\in\bits$,
\begin{equation}\label{eq:P-reversal}
 P_w=P_{\rev w}.
\end{equation}
\end{lem}

\begin{proof}
Word reversal is a bijection
\[
\Sub(a)\longrightarrow\Sub(\rev a),
 \qquad
 u\longmapsto\rev u,
\]
and it preserves both $|u|$ and $b(u)$.  Hence, for every function
$f:\Z\to\mathbb R$, the two right-hand sides in the stopped-tree identity
\eqref{eq:stopped-tree-identity}, applied to $\Sub(a)$ and $\Sub(\rev a)$,
are equal. Thus
\[
 \mathbb E f(S_{\tau_a})=\mathbb E f(S_{\tau_{\rev a}})
\]
for every $f$, which proves \eqref{eq:stopped-reversal}.  Applying
Proposition~\ref{prop:subsequence-representation} first to $w$ and then to
$\rev w$ gives
\[
 P_w=\Law(S_{\tau_{\rev w}})
=\Law(S_{\tau_w})=P_{\rev w}.
\]
\end{proof}

\begin{cor}\label{cor:mu-digit-reversal}
For every positive integer $t$,
\begin{equation}\label{eq:mu-reversal}
\mu_t=\mu_{t^{\mathrm{rev}}},
\end{equation}
where $t^{\mathrm{rev}}$ is obtained by reversing the proper binary expansion of $t$ and ignoring any leading zeros.  For instance, if the binary expansion of $t$ is $110100$, then $t^{\mathrm{rev}}=(001011)_2=(1011)_2$.
\end{cor}

\begin{proof}
If $t=t_w\ge3$ is odd, then $t^{\mathrm{rev}}=t_{\rev w}$, and
Lemma~\ref{lem:convolution} together with \eqref{eq:P-reversal} gives
\[
\mu_t=\mu_1*P_w=\mu_1*P_{\rev w}=\mu_{t^{\mathrm{rev}}}.
\]
The case $t=1$ is immediate. If $t=2^q u$ with $u$ odd, then repeated use of
\eqref{eq:mu-even} gives $\mu_t=\mu_u$, while reversing the $q$ trailing zeros
only creates leading zeros; hence $t^{\mathrm{rev}}=u^{\mathrm{rev}}$.  The
odd case applied to $u$ proves \eqref{eq:mu-reversal}.

This is the binary digit-reversal theorem of Morgenbesser and
Spiegelhofer~\cite{MorgenbesserSpiegelhofer2012}; the argument above gives a
short combinatorial proof within the present subsequence-ideal framework.
\end{proof}

\subsection{Finite stabilization of the nonnegative event}

Under the dictionary $0\leftrightarrow L$ and $1\leftrightarrow R$ of
Remark~\ref{rem:dictionary}, the next proposition is the finite-stabilization
step used in the proof of Lemma~25 of Tar\l{}owski~\cite{Tarlowski2026}.
Indeed, his choice $k=\ell_R(v)+2$, where $\ell_R(v)$ counts the letters
$R$ in $v$, corresponds here to $k=|w|_1+2=s_2(t_w)$, and the equality below
is his stabilization equality in
our notation.  We record it as a separate proposition because its exact form
is used in the final quantitative transfer; the proof follows the same
support-cutoff argument, using the notation developed above.

\begin{prop}\label{prop:event-stabilization}
Let $t=t_w={(1w1)}_2\ge3$ be odd, put
\begin{equation}\label{eq:r-k-definition}
 r=|w|_1,
 \qquad
 k=r+2=s_2(t).
\end{equation}
Then
\begin{equation}\label{eq:event-stabilization}
\mu_t(\Z_{\ge0})=P_{w10^k}(\Z_{\ge0}).
\end{equation}
\end{prop}

\begin{proof}
By \eqref{eq:infinite-comb},
\begin{equation}\label{eq:mu-infinite-sum-event}
\mu_t=\sum_{j=0}^{\infty}2^{-j-1}\sigma_{1-j}P_w.
\end{equation}
At $n=k$, the finite comb formula \eqref{eq:finite-comb} contains the terms
with $0\le j\le k$ and one additional remainder
$2^{-k-1}\sigma_{-k-1}P_w^{(0)}$.

First consider this remainder. Lemma~\ref{lem:support-bound} gives
\[
 \max\supp(P_w^{(0)})\le r+1=k-1.
\]
After translation by $-k-1$, its support is bounded above by
$(k-1)-(k+1)=-2$, so it assigns no mass to $\Z_{\ge0}$.

Next consider an omitted term in \eqref{eq:mu-infinite-sum-event}, for which
$j\ge k+1$. Again by Lemma~\ref{lem:support-bound},
\[
 \max\supp(\sigma_{1-j}P_w)
 \le(r+1)+(1-j)=r+2-j=k-j\le-1.
\]
Thus every omitted term also assigns zero mass to $\Z_{\ge0}$. The finite and infinite measures therefore agree exactly on that event, proving
\eqref{eq:event-stabilization}.
\end{proof}

\begin{proof}[Proof of Theorem~\ref{thm:main}]
Suppose first that $t\ge3$ is odd, and write
\[
 t=t_w={(1w1)}_2.
\]
Let $r$ and $k$ be as in \eqref{eq:r-k-definition}.  By
Proposition~\ref{prop:event-stabilization},
\begin{equation}\label{eq:main-proof-stabilized}
 c_t=P_{w10^k}(\Z_{\ge0}).
\end{equation}
Word-reversal symmetry \eqref{eq:P-reversal} gives
\begin{equation}\label{eq:reverse-final-word}
P_{w10^k}=P_{\rev{(w10^k)}}=P_{0^k1\rev w}.
\end{equation}
By Proposition~\ref{prop:subsequence-representation}, the principal generator
associated with the word $0^k1\rev w$ is its reversal
\begin{equation*}
 \rev{(0^k1\rev w)}=w10^k.
\end{equation*}
This generator ends in $0$. It contains
\begin{equation*}
 O=r+1
 \quad\text{ones},
 \qquad
 Z=|w|_0+k\ge k
 \quad\text{zeros}.
\end{equation*}
Therefore the parameter in Theorem~\ref{thm:strict-principal}(i) is
\begin{equation*}
 m=\min(Z,O+1)
   =\min(|w|_0+k,r+2)
   =k=s_2(t).
\end{equation*}
Combining \eqref{eq:main-proof-stabilized},
\eqref{eq:reverse-final-word}, and
Theorem~\ref{thm:strict-principal}(i), we obtain
\[
 c_t\ge\frac12+2^{-2k-1}
  =\frac12+2^{-2s_2(t)-1}.
\]

For $t=1$, formula \eqref{eq:mu1} gives
\[
c_1=\mu_1(\Z_{\ge0})=\frac12+\frac14=\frac34,
\]
which is stronger than \eqref{eq:main-bound}. Finally, if $t$ is even,
repeated use of \eqref{eq:mu-even} reduces $\mu_t$ to the measure indexed by
the odd part of $t$, while deleting trailing binary zeros does not change
$s_2(t)$. This proves \eqref{eq:main-bound} for every positive integer $t$.
\end{proof}

\section{Further remarks}\label{sec:remarks}

The results operate at three nested levels.  For every finite deletion-closed binary language, Theorem~\ref{thm:ideal-median} places zero at a median of the
first-exit position. For a principal subsequence ideal,
Theorem~\ref{thm:strict-principal} upgrades one side to an explicit strict gap. For the digit-sum correlation trees, Proposition~\ref{prop:event-stabilization}
transfers that finite gap exactly to $c_t$. Thus the subsequence-ideal theorem is not merely a reformulation of Cusick's conjecture: the number-theoretic statement is the final application of a general median principle and its
principal-ideal strengthening.

This separation is conceptually useful. The median theorem belongs to the finite extremal combinatorics of subsequence order ideals, while the final support cutoff belongs to the arithmetic of the sum-of-digits correlation
measures. The proof of Cusick's conjecture needs both parts, but the first
part can be read as a standalone first-exit statement for arbitrary finite deletion-closed languages.

The two finite combinatorial mechanisms are complementary.  The exact formulas
\eqref{eq:exact-plus} and \eqref{eq:exact-minus} reduce the median question to
integer defects at individual balance levels. Marked deletion proves that
each defect is nonnegative, while one forbidden terminal insertion proves
strictness. Separately, word reversal preserves the weighted balance data of a subsequence ideal, yielding the stopped-law symmetry
\eqref{eq:stopped-reversal} and a combinatorial recovery of the classical
digit-reversal theorem.

The finite stabilization is also essential. Total-variation convergence gives convergence of every event probability, but it does not preserve a strict gap above $1/2$.  The support cutoff shows instead that, at the explicit depth
$k=s_2(t)$, both the finite remainder and the entire omitted tail lie strictly below zero. Hence the relevant event probability is already exact at a finite stage.

Natural further questions are to characterize equality in one or both
inequalities in \eqref{eq:ideal-median}, and to identify broader classes of
deletion-closed languages for which a quantitative missing-slot argument
yields strictness. By Proposition~\ref{prop:exact-median-defect}, equality is equivalent to the vanishing of all the relevant level defects. Thus the equality problem is a finite saturation problem for marked insertion slots in
the subsequence order.

\section*{Statements and Declarations}

\noindent\textbf{Competing interests.} The author declares no competing interests.


\begin{thebibliography}{99}

\bibitem{Besineau1972}
J. B\'esineau,
Ind\'ependance statistique d'ensembles li\'es \`a la fonction ``somme des chiffres'',
\emph{Acta Arith.} \textbf{20} (1972), no.~4, 401--416.

\bibitem{CusickLiStanica2011}
T. W. Cusick, Y. Li and P. St\u{a}nic\u{a},
On a combinatorial conjecture,
\emph{Integers} \textbf{11} (2011), no.~2, 185--203.

\bibitem{DrmotaKauersSpiegelhofer2016}
M. Drmota, M. Kauers and L. Spiegelhofer,
On a conjecture of Cusick concerning the sum of digits of $n$ and $n+t$,
\emph{SIAM J. Discrete Math.} \textbf{30} (2016), no.~2, 621--649.

\bibitem{MorgenbesserSpiegelhofer2012}
J. F. Morgenbesser and L. Spiegelhofer,
A reverse order property of correlation measures of the sum-of-digits function,
\emph{Integers} \textbf{12} (2012), Paper No.~A47, 5 pp.

\bibitem{SobolewskiSpiegelhofer2026}
B. Sobolewski and L. Spiegelhofer,
Decomposing the sum-of-digits correlation measure,
\emph{J. Number Theory} \textbf{280} (2026), 702--736.

\bibitem{Spiegelhofer2022}
L. Spiegelhofer,
A lower bound for Cusick's conjecture on the digits of $n+t$,
\emph{Math. Proc. Cambridge Philos. Soc.} \textbf{172} (2022), no.~1, 139--161.

\bibitem{SpiegelhoferWallner2019}
L. Spiegelhofer and M. Wallner,
The Tu--Deng conjecture holds almost surely,
\emph{Electron. J. Combin.} \textbf{26} (2019), no.~1, Paper No.~P1.28, 28 pp.

\bibitem{SpiegelhoferWallner2023}
L. Spiegelhofer and M. Wallner,
The binary digits of $n+t$,
\emph{Ann. Sc. Norm. Super. Pisa Cl. Sci. (5)} \textbf{24} (2023), no.~1, 1--31.

\bibitem{Tarlowski2026}
D. Tar\l{}owski,
On the sum-of-digits measures and Cusick's conjecture via stopped random walks,
arXiv:2605.08624v3, 2026.

\bibitem{TuDeng2011}
Z. Tu and Y. Deng,
A conjecture about binary strings and its applications on constructing
Boolean functions with optimal algebraic immunity,
\emph{Des. Codes Cryptogr.} \textbf{60} (2011), no.~1, 1--14.

\end{thebibliography}
\end{document}